\documentclass[letterpaper, 12pt]{article}

\usepackage[T1]{fontenc}

\usepackage[margin=1in]{geometry}

\usepackage{amscd}

\usepackage{amsmath}

\usepackage{amsfonts}

\usepackage{amssymb}

\usepackage{amsthm}

\usepackage{arydshln}

\usepackage{fancyhdr}

\usepackage{verbatim}

\usepackage{tikz-cd}

\usepackage[all]{xy}

\usepackage{color}

\usepackage[colorlinks=true, linkcolor=blue, citecolor=blue,
pagebackref=true]{hyperref}

\usepackage{fouriernc}

\usepackage{mathtools} \usepackage{etoolbox}
\theoremstyle{definition}

\theoremstyle{plain} \newtheorem{theorem}{Theorem}[section]

\theoremstyle{plain} 

\theoremstyle{plain} 

\theoremstyle{plain} 

\theoremstyle{plain} 

\theoremstyle{remark} \newtheorem*{remark}{Remark}

\theoremstyle{definition} 

\theoremstyle{definition} \newtheorem*{definition*}{Definition}

\theoremstyle{definition} \newtheorem{question}[theorem]{Question}

\theoremstyle{remark} 


\newcommand{\PP}{\mathbb{P}}

\newcommand{\RR}{\mathbb{R}}

\newcommand{\NN}{\mathbb{N}}

\newcommand{\TT}{\mathbb{T}}

\newcommand{\ZZ}{\mathbb{Z}}

\newcommand{\calA}{\mathcal{A}}

\newcommand{\bone}{\mathbf{1}}

\renewcommand{\leq}{\leqslant} \renewcommand{\geq}{\geqslant}

\DeclarePairedDelimiter{\abs}{\lvert}{\rvert}
\DeclarePairedDelimiter{\norm}{\lVert}{\rVert}

\DeclarePairedDelimiter{\set}{\lbrace}{\rbrace}
\DeclarePairedDelimiter{\parens}{\lparen}{\rparen}

\DeclareMathOperator{\meas}{meas}

\def\eps{{\varepsilon}}

\def\1int{{[0,1]}}

\title{A note on sequences not having metric Poissonian pair correlations}

\author{Felipe A.~Ram{\'i}rez\footnote{\texttt{framirez@wesleyan.edu}}
  \\ Wesleyan University}

\date{}

\begin{document}

\maketitle

% ===============================================
\begin{abstract}
  The purpose of this note is to present a construction of sequences
  which do not have metric Poissonian pair correlations (MPPC) and
  whose additive energies grow at rates that come arbitrarily close to
  a threshold below which it is believed that \emph{all} sequences
  have MPPC. A similar result appears in work of Lachmann and Technau
  and is proved using a totally different strategy. The main novelty
  here is the simplicity of the proof, which we arrive at by modifying
  a construction of Bourgain.
\end{abstract}
% ===============================================

\thispagestyle{empty}

% \setcounter{tocdepth}{1} %excludes subsections from toc
% %\tableofcontents
% {\footnotesize{\tableofcontents}}

\paragraph*{Notation}

Let us set some notation once and for all. For functions $f,g : \NN\to \RR_{\geq 0}$, we write $f(n)\ll g(n)$ if there is a constant $c>0$ such that $f(n) \leq c g(n)$ holds for all sufficiently large $n\in\NN$. We write $f(n) \asymp g(n)$ if $f(n) \ll g(n)$ and $g(n) \ll f(n)$ both hold. We write $f(n)\sim g(n)$ if $\lim_{n\to\infty} f(n)/g(n) = 1$. 

\section{Introduction}
\label{sec:introduction}

Let $\calA\subset \NN$ be an infinite subset and denote its smallest
$N$ elements $A_N$. For $N\in\NN$, $\alpha\in[0,1]$, and $s>0$, the
quantity
\begin{equation}\label{eq:2}
  F(\alpha, s, N, \calA) = \frac{1}{N}
  \sum_{\substack{(a,b)\in A^2 \\ a\neq b}} \bone_{[-s/N,s/N] + \ZZ}(\alpha(a-b)),
\end{equation}
measures how often two points in
$\alpha A_N (\bmod 1)$ lie within a distance $2s/N$ of each other on
the circle $\TT = \RR/\ZZ$. Grouping the terms in the sum according to the differences $a-b$ leads to the equivalent and convenient expression
\begin{equation}\label{eq:alternate}
  F(\alpha, s, N, \calA) = \frac{1}{N}
  \sum_{d\in\ZZ\setminus\set{0}} \abs*{A_N\cap (A_N + d)}\,\bone_{d, s/N}(\alpha),
\end{equation}
where $\bone_{d,\eps}$ denotes the indicator function of the set
$\set{\alpha\in[0,1] : \norm{d\alpha}\leq \eps}$, where $\norm{\cdot}$ is distance to $\ZZ$ . 

If for almost every $\alpha\in[0,1]$ we
have $F(\alpha, s, N, \calA)\sim 2s$, then $\calA$ is said to have
\emph{metric Poissonian pair correlations (MPPC)}. Since a random
point sequence on $\TT$ will almost surely have asymptotically
Poissonian pair correlations, MPPC is understood as a property
connoting random-like behavior for an integer sequence. It is of great
interest to understand which integer sequences do and do not have
metric Poissonian pair correlations.

In~\cite{RudnickSarnak}, Rudnick and Sarnak showed that the sequence
$\parens{n^k}_{n=1}^\infty$ has MPPC whenever $k\geq 2$, whereas it is
easy to show that it does not have MPPC if $k=1$. The intuitive reason
that $\parens{n}_{n=1}^\infty$ does not have MPPC is that in this case
$A_N = \set{1, \dots, N}$, and one quickly sees that the quantities
$(a-b)$ arising in~\eqref{eq:2} are too structured to be
random-like. Aistleitner, Larcher, and Lewko made this intuition
rigorous by connecting MPPC to the behavior of
\begin{equation*}
  E(A_N) = \#\set*{(a,b,c,d)\in A_N^4 : a +b = c + d},
\end{equation*}
the \emph{additive energy} of $A_N$. They showed that $\calA$ has MPPC
whenever there exists some $\eps>0$ for which $E(A_N) \ll N^{3-\eps}$
holds~\cite{Aistleitneretaladditiveenergy}. In the appendix to the
same paper, Bourgain showed that if $E(A_N)\gg N^3$ then $\calA$ does
not have MPPC, and also that there exist sequences for which
$E(A_N) = o(N^3)$ which do not have MPPC. This led to a series of
papers exploring the connection between additive energy and
MPPC. Bloom, Chow, Gafni, and Walker asked the
following guiding question.

\begin{question}[{\cite[Fundamental Question 1.7]{Bloometal}}]\label{q}
  Suppose there is a nonincreasing $\psi:\NN\to [0,1]$ such that
  $E(A_N) \sim N^3\psi(N)$. Is convergence of $\sum \psi(N)/N$
  necessary and sufficient for the sequence $\calA$ to have metric
  Poissonian pair correlations?
\end{question}

They proved results in support of the answer being ``yes,'' but as of
this writing the overall picture is not complete. 

For the sufficiency part, the best result so far is due to Bloom and
Walker, and it says that there exists some universal constant $C>1$
such that if $\psi(N) \ll (\log N)^{-C}$, then $\calA$ has
MPPC~\cite{BloomWalker}. Of course, if one believes the sufficiency
part of Question~\ref{q}, then one should believe that any $C>1$ will
do. (Indeed, Hinrichs \emph{et al.}~have established this for a higher
dimensional version of the problem~\cite{HKLSUPPC}.)

The answer to the necessity part of Question~\ref{q} turns out to be
``no.'' Aistleitner, Lachmann, and Technau found, for any $\eps>0$, sequences
$\calA\subset\NN$ for which 
\begin{equation*}
  E(A_N) \gg N^3 (\log N)^{-\frac{3}{4} - \eps},
\end{equation*}
yet they have metric Poissonian pair correlations~\cite{AistleitnerLachmannTechnau}. However, the
construction is very special. There is still reason to think that
perhaps a ``randomly chosen'' sequence $\calA\subset\NN$ whose
additive energy behaves as in the divergence part of Question~\ref{q}
will not have MPPC. Bloom \emph{et al.}~proved a result to this
effect, showing that in a certain random model, a sequence $\calA$
whose additive energy satisfies
\begin{equation*}
  E(A_N) \asymp N^3 (\log N)^{-1}(\log\log N)^{-C}
\end{equation*}
for some $0\leq C\leq 1$ will almost surely \emph{not} have
MPPC~\cite[Theorem~1.6]{Bloometal}. In~\cite[Theorem 2]{LachmannTechnau},
Lachmann and Technau constructed examples where the additive energy is of
order $E(A_N)\asymp N^3\psi(N)$ where $\psi$ is any function as in the
divergence part of Question~\ref{q} that satisfies the further
condition that $\psi(N) \gg N^{-1/3}(\log
N)^{7/3}$. In particular, this yields examples of sets
$\calA\subset\NN$ where
\begin{equation*}
E(A_N) \asymp N^3(\log N \log\log N \dots
\underbrace{\log\log\dots\log}_{r \textrm{ iterates}} N)^{-1}
\end{equation*}
which do not have metric Poissonian pair correlations. 

\

In this note, we present a modified version of Bourgain's
construction~\cite[Appendix]{Aistleitneretaladditiveenergy} which
gives examples of sequences which do not have MPPC and whose additive
energies meet the threshold proposed in Question~\ref{q}. That is, we
prove the following.

\begin{theorem}\label{thm}
  Suppose $\psi:\NN\to[0,1]$ is a nonincreasing function such that
  $N^{3-\delta}\psi(N)$ is nondecreasing for some fixed $\delta>0$,
  and such that $\sum\psi(N)/N$ diverges. Then there exists an
  infinite set $\calA\subset\NN$ such that $E(A_N)\asymp N^3\psi(N)$
  and such that $\calA$ does not have metric Poissonian pair
  correlations.
\end{theorem}

\begin{remark}
  As in~\cite[Theorem 2]{LachmannTechnau}, Theorem~\ref{thm} has a
  condition on $\psi$ besides just divergence of the series. Since
  $E(A_N)$ must increase to infinity it is unavoidable that such a
  theorem should have extra conditions on $\psi$. Indeed, the extra
  condition in Theorem~\ref{thm} is only used in the proof that the
  constructed sequence $\calA$ actually satisfies
  $E(A_N)\ll N^3\psi(N)$. It is not used in proving the assertion that
  $\calA$ does not have MPPC.

  Given that there has to be some extra condition on $\psi$, perhaps it would be
  most natural to only require that $N^3\psi(N)$ increase to
  infinity. Instead, we make the slightly stronger assumption that
  there is some $\delta>0$ for which $N^{3-\delta}\psi(N)$ is
  nondecreasing. This is also not so unnatural, since the divergence of
  $\sum\psi(N)/N$ already requires that $N^{3-\delta} \psi(N)$ is
  unbounded whenever $0 < \delta < 3$. In particular,
  Theorem~\ref{thm} applies when
\begin{equation*}
\psi(N) =(\log N \log\log N \dots
\underbrace{\log\log\dots\log}_{r \textrm{ iterates}} N)^{-1}
\end{equation*}
as in~\cite{LachmannTechnau}. 
\end{remark}

The rest of this note consists of the proof. For a more detailed
discussion of pair correlations and additive energy, we recommend the
surveys~\cite{LarcherStockinger, WalkerSurvey}.

\section{Proof of Theorem~\ref{thm}}
\label{sec:proof}

Notice that we lose no generality in assuming that $\psi(N)=o(1)$, for
otherwise we would have $E(A_N)=\Omega(N^3)$, and in this case it is
known that $\calA$ cannot have MPPC. We may also assume that
$\psi(N)^{-1}$ takes only integer values. 

Let $\iota(N) :\NN\to \RR$ decrease to $0$ slowly enough that
$\sum\frac{\psi(N)\iota(N)}{N}$ still diverges. Let $(\Delta_N)_N$ be a
positive integer sequence that increases fast enough that the sets
\begin{equation*}
  S_N :=\set*{\alpha\in[0,1] : \norm{\Delta_N d \alpha} \leq
    \frac{\psi(N)\sqrt{\iota(N)}}{N} \quad\textrm{ for some }\quad 0 < d \leq N\sqrt{\iota(N)}}
\end{equation*}
are pairwise quasi-independent, meaning that there is some constant $C>0$ such that 
\begin{equation*}
\meas(S_N\cap S_M) \leq C\meas(S_N)\meas(S_M) \qquad (M\neq N).
\end{equation*}
To see that it is possible to do this,
note that $S_N = \Delta_N^{-1} S$ where $S$ is a union of finitely
many intervals in $\TT$. In particular, $S$ is measurable. Recall that
for any $m\geq 2$, the ``times $m$ modulo $1$'' map $T_m:\TT\to\TT$ is
measure-preserving, meaning that for any measurable set $S$ we have
$\meas(T_m^{-1}S)=\meas(S)$, and mixing, meaning that for any two measurable sets
$S, T\subset\TT$ we have
\begin{equation*}
  \lim_{k\to\infty} \meas\parens*{T_m^{-k}(S)\cap T} = \meas(S)\meas(T).
\end{equation*}
We may therefore take $\Delta_1=1$ and inductively set $\Delta_N$ to be a power of $m$ that is large enough that
\begin{equation*}
\meas(S_N\cap S_M) \leq 2\meas(S_N)\meas(S_M)
\end{equation*}
for all $M<N$.

Notice that $\meas(S_N)\gg \psi(N)\iota(N)$. Since
$\sum_N\frac{\psi(N)\iota(N)}{N}$ is a divergent sum of nonincreasing
terms, by Cauchy's condensation test we have that
$\sum_t \psi(2^t)\iota(2^t)$ diverges, hence $\sum_t \meas(S_{2^t})$
diverges. Since the sets $(S_{2^t})_t$ are pairwise quasi-independent,
the version of the second Borel--Cantelli lemma proved by Erd{\H
  o}s--Renyi~\cite{ErdosRenyiBC} guarantees that the limsup
set $S_\infty:=\limsup_{t\to\infty}S_{2^t}$ has positive measure.

Our goal now is to construct a sequence $\calA\subset\NN$ such that
$E(A_N)\asymp N^3\psi(N)$ and such that for every $\alpha\in S_\infty$, we have
$\limsup_{N\to\infty}F(\alpha, 1, N, \calA)=\infty$. We will construct $\calA$ block by block. For each $N$, let
\begin{equation*}
  B_N(\omega) = \set*{\Delta_N\parens*{\frac{N}{\psi(N)} + n} : 1 \leq n \leq
    \frac{N}{\psi(N)}\quad\textrm{and}\quad \xi_n^{(N)}(\omega)=1},
\end{equation*}
with $\xi_1^{(N)}, \dots, \xi_{N/\psi(N)}^{(N)}$ independent Bernoulli
random variables with $\PP(\xi_n^{(N)}=1)=\psi(N)$. (Recall that we have assumed without loss of generality that $\psi(N)^{-1}$ is always an integer, so the blocks $B_N(\omega)$ consist only of integers.) For comparison, these
blocks $B_N(\omega)$ are dilates of the blocks
in~\cite{Aistleitneretaladditiveenergy} by the factor $\Delta_N$.

In light of~\cite[Lemma 6]{Aistleitneretaladditiveenergy}, the
following three properties hold with positive probability, and so we
may henceforth assume that $B_N$ is an instantiation of $B_N(\omega)$
where:
\begin{enumerate}
\item For all $d\in \ZZ\setminus\set{0}$ we have $\abs{B_N\cap (B_N +
    \Delta_Nd)}\leq 2N\psi(N)$. 
\item For all $d\in \ZZ\setminus\set{0}$ with $\abs{d} <
  \frac{N}{10\psi(N)}$ we have $\abs{B_N\cap (B_N +
    \Delta_Nd)}\geq \frac{1}{2}N\psi(N)$.
\item We have $N/2 \leq \abs{B_N} \leq 2N$. 
\end{enumerate}
Since any two elements of $B_N$ differ by a multiple of $\Delta_N$, we
have
\begin{equation*}
  E(B_N)  = \sum_{d\in\ZZ}\abs*{B_N\cap (B_N + \Delta_N d)}^2.
\end{equation*}
With this, the first two properties above show us that
$E(B_N) \asymp N^3\psi(N)$.

Let $\calA :=\set{B_1, B_2, B_4, \dots, }$ be the concatenation of the
blocks $B_{2^t}$, $t\geq 0$. Suppose that $A_N$ is a truncation of
$\calA$ in the block $B_{2^t}$. It is obvious then that
\begin{equation*}
  E(A_N) \geq E(B_{2^{t-1}}) \gg (2^{t-1})^3 \psi(2^{t-1}) \gg N^3
  \psi(N). 
\end{equation*}
Also, by possibly making $(\Delta_N)_N$ sparser if needed, we have
\begin{align*}
  E(A_N) &\leq \sum_{k=0}^t E(B_{2^k}) \\
         &\ll \sum_{k=0}^t 2^{3k}\psi(2^k)\\
         &\ll 2^{3t}\psi(2^t) \\
  &\ll N^3\psi(N),
\end{align*}
where we have used our assumption that $N^{3-\delta}\psi(N)$ is nondecreasing in the third line. This shows that $\calA$ has the desired behavior in additive energy,
namely, $E(A_N)\asymp N^3\psi(N)$.

As for pair correlations, recall from~\eqref{eq:alternate} that
\begin{equation*}
  F(\alpha, s, N, \calA) = \frac{1}{N}
  \sum_{d\in\ZZ\setminus\set{0}} \abs*{A_N\cap (A_N + d)}\,\bone_{d, s/N}(\alpha).
\end{equation*}
In particular, for $N=2^t$ we have
\begin{align}
  F(1, \abs{B_1}+ \abs{B_2} + \abs{B_4}+\dots+\abs{B_N}, \calA) &\geq
                                                                  \frac{1}{4N} \sum_{d\neq 0} \abs{B_N \cap (B_N
                                                                  +\Delta_Nd)}\bone_{\Delta_Nd,1/(4N)} \nonumber \\
                                                                &\geq
                                                                  \frac{\psi(N)}{8}\sum_{0
                                                                  <
                                                                  \abs{d}\leq
                                                                  \frac{N}{10\psi(N)}}
                                                                  \bone_{\Delta_N
                                                                  d,1/(4N)}
\nonumber  \\
                                                                &\geq
                                                                  \frac{\psi(N)}{4}\sum_{0
                                                                  <
                                                                  d\leq
                                                                  \frac{N}{10\psi(N)}}
                                                                  \bone_{\Delta_N
                                                                  d,1/(4N)}. \label{eq:sum}
\end{align}
Notice that for any $\alpha\in S_N$ there is some $0 < d \leq N\sqrt{\iota (N)}$ for which $\norm{\Delta_N d \alpha} \leq \frac{\psi(N)\sqrt{\iota(N)}}{N}$. Each of the positive multiplies $d, 2d, 3d, \dots, kd$ will contribute $1$ to the sum in~\eqref{eq:sum} as long as $kd \leq \frac{N}{10\psi(N)}$ and $k\norm{\Delta_N d \alpha} \leq \frac{1}{4N}$. We are therefore assured at least $\frac{1}{10\psi(N)\sqrt{\iota(N)}}$ contributions, hence
\begin{equation*}
  F(\alpha, 1, \abs{B_1}+\dots+\abs{B_N}, \calA) \geq
  \frac{1}{40\sqrt{\iota(N)}}
\end{equation*}
for every $\alpha\in S_N$. Since a positive-measure set of $\alpha \in [0,1]$ is contained in infinitely many
$S_{2^t}$'s, this implies that 
\begin{equation*}
\limsup_{N\to\infty} F(\alpha, 1, N, \calA) = \infty
\end{equation*}
for a positive-measure set of $\alpha\in [0,1]$.  Therefore $\calA$ does not have
MPPC.\qed

\subsection*{Acknowledgments}
\label{sec:acknowledgments}

I thank Niclas Technau for an informative conversation and the anonymous referee for making suggestions that improved the presentation.

\bibliographystyle{plain}

%\bibliography{../bibliography}

\end{document}